\documentclass[journal]{IEEEtran}      

\pdfoutput=1
\IEEEoverridecommandlockouts                              

\newtheorem{Remark}{Remark}

\newtheorem{Definition}{Definition}

\newenvironment{Proof}{\noindent{\em Proof:\/}}{\hfill $\Box$\par}
\newtheorem{Theorem}{Theorem}
\newtheorem{Lemma}{Lemma}

\newtheorem{Assumption}{Assumption}

\usepackage{graphicx}
\usepackage{amsmath,bm}
\usepackage{amssymb}
\usepackage{algorithm}
\usepackage{algpseudocode}
\usepackage{cite}
\usepackage{hyperref}
\usepackage[caption=false]{subfig}
\hypersetup{hidelinks}

\algdef{SE}[DOWHILE]{Do}{doWhile}{\algorithmicdo}[1]{\algorithmicwhile\ #1}%

\newcommand{\mathactivatecomma}{%
  \begingroup\lccode`~=`\,
  \lowercase{\endgroup\edef~}{\mathchar\the\mathcode`\,\penalty0 }}

\algnewcommand{\Initialize}[1]{%
  \State \textbf{Initialize: $i \in \mathcal{V}$}
  \Statex \hspace*{\algorithmicindent}\parbox[t]{.8\linewidth}{\raggedright #1}
}

\algnewcommand{\Iteration}[1]{%
  \State \textbf{Iteration $(k\geq 0)$: $i \in \mathcal{V}$}
  \Statex \hspace*{\algorithmicindent}\parbox[t]{.8\linewidth}{\raggedright #1}
}

\algnewcommand{\Output}[1]{%
  \State \textbf{Output: $i \in \mathcal{V}$}
  \Statex \hspace*{\algorithmicindent}\parbox[t]{.8\linewidth}{\raggedright #1}
}

\title{\LARGE \bf
Distributed Nash Equilibrium Seeking with Limited Cost Function Knowledge via A Consensus-Based Gradient-Free Method*
}

\author{Yipeng Pang and Guoqiang Hu
\thanks{*This research was supported in part by Singapore Ministry of Education Academic Research Fund Tier 1 RG180/17(2017-T1-002-158), and in part by the National Research Foundation, Prime Minister's Office, Singapore under the Energy Innovation Research Programme (EIRP) for Building Energy Efficiency Grant Call, administered by the Building and Construction Authority (NRF2013EWT-EIRP004-051).}
\thanks{Y. Pang and G. Hu are with the School of Electrical and Electronic Engineering, Nanyang
Technological University, 639798, Singapore
        {\tt\small ypang005@e.ntu.edu.sg, gqhu@ntu.edu.sg}.}%
}

\begin{document}

\bstctlcite{IEEEexample:BSTcontrol}

\maketitle
\thispagestyle{empty}
\pagestyle{empty}

\begin{abstract}
This paper considers a distributed Nash equilibrium seeking problem, where the players only have partial access to other players' actions, such as their neighbors' actions. Thus, the players are supposed to communicate with each other to estimate other players' actions. To solve the problem, a leader-following consensus gradient-free distributed Nash equilibrium seeking algorithm is proposed. This algorithm utilizes only the measurements of the player's local cost function without the knowledge of its explicit expression or the requirement on its smoothness. Hence, the algorithm is gradient-free during the entire updating process. Moreover, the analysis on the convergence of the Nash equilibrium is studied for the algorithm with both diminishing and constant step-sizes, respectively. Specifically, in the case of diminishing step-size, it is shown that the players' actions converge to the Nash equilibrium almost surely, while in the case of fixed step-size, the convergence to the neighborhood of the Nash equilibrium is achieved. The performance of the proposed algorithm is verified through numerical simulations.
\end{abstract}


\begin{IEEEkeywords}
Nash equilibrium seeking, non-cooperative games, game theory, gradient-free methods, distributed algorithms.
\end{IEEEkeywords}
\section{Introduction}
Over the decades, game theory, as a power tool of analyzing the strategic interactions between rational decision-makers, has found its great potential in various application fields such as social science, economics, electricity markets, power systems, to list a few. An important concept in game theory, Nash equilbrium, named after John Forbes Nash Jr., is a proposed solution in non-cooperative games involving two or more players. Recently, with the emergence of multi-agent system, Nash equilibrium seeking in multi-player non-cooperative games has received increasing attention. More precisely, this type of games involves a number of players, who selfishly minimize their own cost functions by making decisions in response to other players' actions. 

Recently, a large number of studies on Nash equilibrium computation in non-cooperative games have been reported, 
such as \cite{Lou2016,Tatarenko2018,Romano2018,Belgioioso2018,Deng2019a,Yi2019} to list a few. 
The challenge of such problem settings is the requirement of global knowledge on all players' actions, which is not practical if the underlying communication network is not fully connected. In such cases, players have to make decisions based only on a limited set of information, such as the information from the neighbors. Therefore, a distributed information sharing protocol is usually adopted to disseminate the local information among players. For example, a dynamic average consensus protocol was adopted in \cite{Ye2017a}, where a primal-dual dynamic based seeking strategy was developed to find Nash equilibrium in set constrained aggregate games. It was also utilized in \cite{Deng2019} with the help of differential inclusions and differentiated projections for aggregative games, where the players' actions are coupled by linear constraints. The dynamic average consensus protocol was also proposed to achieve simultaneous social cost minimization and Nash equilibrium in a class of $N$-coalition games in \cite{Ye2018}. Different from these works, the work in \cite{Sun2018} considered a continuous time generalized convex game with shared inequality constraints among players, and proposed a leader-following consensus protocol with gradient descent method to compute the generalized Nash equilibrium. This protocol was also employed in \cite{Lu2019} to estimate the other players' actions for the generalized games, where the players' action sets are constrained by nonlinear inequality and linear equations. Apart from the leader-following consensus and dynamic average consensus protocols, gossip-based averaging techniques were also commonly utilized in Nash equilibrium computations, such as \cite{Koshal2012,Salehisadaghiani2016}.
Most of the existing literature including the aforementioned works are model-based approaches, \textit{i.e.}, the implementation of the algorithms relies on the knowledge of the explicit form of the players' cost functions, such as the derivative computation. However, the requirement of the knowledge on the explicit expression of players' cost functions is restrictive in the cases where the input/output relationship is difficult to model.

There are non-model based approaches, which utilize the players' local measurements without the requirement on the information of the functional form. For example, the work in \cite{Zhu2016} considered a generalized convex game with both convex coupling inequality constraints and local set constraints. A finite-differencing method with two-way perturbations was proposed to approximate the partial gradient. The perturbation parameter needs to be chosen carefully to match the selected step-size. Different from that, the work in \cite{Tatarenko2019} proposed a distributed payoff-based algorithm for a class of convex games with and without coupling constraints. This technique was further extended in \cite{Tatarenko2019a} where the algorithm convergence was proved under mere monotonicity assumption. Overall, the payoff-based learning strategy proposed in these two works enables players to sample their actions in a Gaussian distribution. Then, the mean of this distribution is iteratively updated using only local payoff values. Another typical non-model based approaches are extremum seeking-based methods, such as \cite{Liu2011,Stankovic2012,Frihauf2012,Ye2014,Ye2016a,Ye2015a}. Specifically, the work in \cite{Liu2011} proposed a continuous time multi-input stochastic extremum seeking algorithm for the Nash equilibrium seeking in non-cooperative games with general nonlinear cost functions. In \cite{Stankovic2012}, a discrete time stochastic extremum seeking method was presented in non-cooperative games where the players' cost functions are strictly convex, but the actions are subject to a linear dynamic constraint. The work in \cite{Frihauf2012} developed an integrator-type extremum seeking algorithm in non-cooperative games with both quadratic payoffs and general non-quadratic payoffs as the output of a dynamic system. More extremum seeking algorithms have been proposed in potential games with unstable dynamics \cite{Ye2014}, dynamical constraints \cite{Ye2016a}, and non-cooperative games with time-varying Nash equilibrium \cite{Ye2015}. In general, the extremum seeking strategy makes use of the cost value together with some sinusoidal dither signals for perturbation, such that the gradient of the cost function is extractable. Even though all the aforementioned works need no explicit model information during the implementation, they assume the players' cost functions to be smooth to some extent, which can be restrictive if the players' cost functions are generally non-differentiable. This motivates the study of gradient-free technique, which is free of the knowledge on the explicit expressions of the players' cost functions and applicable to non-differentiable problems. In fact, gradient-free algorithms have been studied in distributed optimization problems \cite{Li2015,Yuan2015,Chen2017,Pang2017,Pang2018,Pang2020}. However, little attention has been received in non-cooperative games.

In this paper, we focus on the research of non-model based Nash equilibrium seeking methods. Specifically, a gradient-free distributed algorithm is proposed to solve the Nash equilibrium seeking problem in a multi-player non-cooperative static game under a directed communication graph. As compared to the existing literature, the major contributions of this paper are twofold.


\begin{enumerate}
\item The proposed algorithm does not rely on the knowledge of the explicit form of the players' cost functions. Different from non-model based approaches such as payoff-based learning \cite{Tatarenko2019,Tatarenko2019a} and extremum seeking \cite{Liu2011,Stankovic2012,Frihauf2012,Ye2014,Ye2016a,Ye2015a}, the proposed algorithm allows the cost functions to be non-smooth. Unlike the finite-differencing method in \cite{Zhu2016} where the perturbation parameter needs to match the step-size, the proposed algorithm establishes the convergence to the Nash equilibrium with only the requirement of a small smoothing parameter.
\item The convergence of the proposed algorithm to the Nash equilibrium is rigorously studied for both diminishing and constant step-sizes, respectively. Specifically, for the diminishing case, an exact convergence to the Nash equilibrium is attained, while for the constant case, an approximate convergence to the Nash equilibrium with the gap proportional to the step-size is achieved.
\end{enumerate}

The paper is organized as follows. The problem is defined in section~\ref{sec:problem_formulation}. Main procedures of the proposed algorithm are described in section~\ref{sec:distr_opt}. The convergence analysis of the proposed algorithm for both diminishing step-size and constant step-size is presented in section~\ref{sec:conv_analysis}. In section~\ref{sec:simulation}, the performance of the proposed algorithm is illustrated through a numerical example. Section~\ref{sec:conclusion} concludes the paper.

\section{Problem Formulation}\label{sec:problem_formulation}

This section firstly introduces the notations used. Then, the problem is formulated, followed by some preliminary results.

\subsection{Notations}
We use $\mathbb{R}$ and $\mathbb{R}^N$ to denote the set of real numbers and $N$-dimensional column vectors, respectively. 
For a matrix $A$, the element in the $i$-th row and $j$-th column of $A$ is represented by $[A]_{ij}$, and its transpose is denoted by $A^\top$. We write $\mathbb{E}[\cdot]$ to denote the expected value of a random variable. 
For any two vectors $\mathbf{x}$ and $\mathbf{y}$, the operator $\langle\mathbf{x},\mathbf{y}\rangle$ denotes the inner product of $\mathbf{x}$ and $\mathbf{y}$. We use $\|\mathbf{x}\|$ for the standard Euclidean norm of a vector $\mathbf{x}$, \textit{i.e.}, $\|\mathbf{x}\| = \sqrt{\langle\mathbf{x},\mathbf{x}\rangle}$, and $\mathcal{P}_{\Omega}[\mathbf{x}]$ for the projection of a vector $\mathbf{x}$ on the set $\Omega$, \textit{i.e.}, $\mathcal{P}_{\Omega}[\mathbf{x}] = \arg\min_{\hat{\mathbf{x}} \in \Omega} \|\hat{\mathbf{x}} - \mathbf{x}\|^2$.
For a differentiable function $f$, we use $\nabla_x f(x,y)$ to represent its partial derivative with respect to $x$ at the point $(x,y)$.
For a possibly non-differentiable function $f$, we denote its $\varepsilon$-subdifferential by $\partial^\varepsilon_{x}f(x,y)$ at $x$ for any fixed $y$, \textit{i.e.}, $f(z,y) \geq f(x,y) - \varepsilon + \langle g, z-x \rangle$, $g \in \partial^\varepsilon_{x}f(x,y)$. If $\varepsilon=0$, we simplify the notation to $\partial_{x}f(x,y)$, which is the set of its subgradients at $x$ for any fixed $y$.
For a sequence of random vectors $a_k$, we say that $a_k$ converges to $a$ almost surely, if the probability of $\lim_{k\to\infty}a_k=a$ is 1.

\subsection{Problem Definition}
We consider a directed communication graph represented by $\mathcal{G} = \{\mathcal{V},\mathcal{E}\}$, where $\mathcal{V}$ is the set of agents, and $\mathcal{E} \subset \mathcal{V}\times\mathcal{V}$ is the set of edges, \textit{i.e.}, for any $i,j\in\mathcal{V}$, the ordered pair $(i,j)\in\mathcal{E}$ if and only if the information can be transfered from agent $i$ to agent $j$. In particular, the set $\mathcal{E}$ includes $(i,i)$ for all $i\in\mathcal{V}$. A matrix $A$ associated with the directed graph $\mathcal{G}$ is known as the adjacency matrix, which is designed such that $[A]_{ij}> 0$ if $(j,i)\in\mathcal{E}$ and $[A]_{ij}= 0$ otherwise. The set of in-neighbors (respectively, out-neighbors) of agent $i$ is denoted by $\mathcal{N}^{\text{in}}_i = \{j \in \mathcal{V} | (j,i)\in \mathcal{E}\}$ (respectively, $\mathcal{N}^{\text{out}}_i = \{j \in \mathcal{V} | (i,j)\in \mathcal{E}\}$). In particular, agent $i$ is both an in-neighbor and an out-neighbor of itself, \textit{i.e.}, $i\in\mathcal{N}^{\text{in}}_i$ and $i\in\mathcal{N}^{\text{out}}_i$. It should be noted that $\mathcal{N}^{\text{in}}_i\neq\mathcal{N}^{\text{out}}_i$ in general. 

Consider a game $\Gamma$ with $N$ players that communicate with each other under a directed communication graph $\mathcal{G}= \{\mathcal{V},\mathcal{E}\}$. The set of players is $\mathcal{V} = \{1, 2, \ldots, N\}$. Each player $i\in \mathcal{V}$ owns a cost function $f_i:\Omega\to\mathbb{R}$, where $\Omega =  \Omega_i\times\Omega_{-i}\subset \mathbb{R}^N$ is the action set of all players, and $\Omega_i\subset\mathbb{R}$ (respectively, $\Omega_{-i}\subset\mathbb{R}^{N-1}$) denotes the action set of player $i$ (respectively, all players except player $i$). Let $\mathbf{x} = (x_i,\mathbf{x}_{-i})\in\Omega$ be the vector of all players' actions, where $x_i\in\Omega_i$ (respectively, $\mathbf{x}_{-i}\in\Omega_{-i}$) represents the action of player $i$ (respectively, all players except player $i$). Under the communication graph $\mathcal{G}$, if player $j$ is not an out-neighbor of player $i$ (\textit{i.e.}, $j\notin\mathcal{N}^{\text{out}}_i$), then player $j$ does not have direct access to player $i$'s action. Game $\Gamma(N,\{f_i\},\{\Omega_i\},\mathcal{G})$ is played such that for given $\mathbf{x}_{-i}\in\Omega_{-i}$, the objective of each player $i\in \mathcal{V}$ is to minimize its own cost function, \textit{i.e.},
\begin{align}\label{eq:cost_function}
  \min_{x_i\in\Omega_i} f_i(x_i,\mathbf{x}_{-i}), \quad i\in \mathcal{V}.
\end{align}
It should be emphasized that the explicit mathematical expression of the cost function $f_i$ is unknown, but each player $i\in\mathcal{V}$ can measure the value of $f_i$ by introducing some input to the system $f_i$. Moreover, the solution set of player $i\in\mathcal{V}$ to the problem \eqref{eq:cost_function} is dependent on the other players' action $\mathbf{x}_{-i}$, which may not be directly accessible. Thus, the objective is to develop a distributed strategy such that all players' actions converge to a Nash equilibrium under the communication graph $\mathcal{G}$.

The formal definition of a Nash equilibrium of a game $\Gamma$ is given below \cite{Nash1951}.
\begin{Definition}\label{definition_NE}
Consider a game $\Gamma(N,\{f_i\},\{\Omega_i\},\mathcal{G})$. Nash equilibrium is an action profile where no player can reduce its cost by unilaterally changing its own action, \textit{i.e.}, a vector $\mathbf{x}^\star = (x_i^\star,\mathbf{x}_{-i}^\star)\in\Omega$ is called a Nash equilibrium of the game $\Gamma(N,\{f_i\},\{\Omega_i\},\mathcal{G})$ if and only if
\begin{align*}
f_i(x_i^\star,\mathbf{x}_{-i}^\star)\leq f_i(x_i,\mathbf{x}_{-i}^\star),\quad \forall x_i\in\Omega_i, \forall i\in \mathcal{V}.
\end{align*}
\end{Definition}

The following standard assumptions are made throughout the paper.
\begin{Assumption}\label{assumption_graph}
The directed graph $\mathcal{G}$ is strongly connected and its associated adjacency matrix $A$ is doubly-stochastic, \textit{i.e.}, $\sum_{j=1}^N[A]_{ij} = 1$ for all $i \in \mathcal{V}$, and $\sum_{i=1}^N[A]_{ij} = 1$ for all $j \in \mathcal{V}$.
\end{Assumption}

\begin{Assumption}\label{assumption_local_f_lipschitz}
For each player $i\in\mathcal{V}$, its action set $\Omega_i$ is non-empty, convex and compact. The cost function $f_i(x_i,\mathbf{x}_{-i})$ is convex in $x_i$ for every $\mathbf{x}_{-i}$, and jointly continuous in $\mathbf{x}$ but not necessarily differentiable. Also, $f_i(x_i,\mathbf{x}_{-i})$ is Lipschitz continuous in $x_i$ (respectively, $\mathbf{x}_{-i}$) for every fixed $\mathbf{x}_{-i}$ (respectively, $x_i$),
\textit{i.e.}, $\forall x_i, y_i \in \Omega_i$ (respectively, $\forall \mathbf{x}_{-i}, \mathbf{y}_{-i} \in \Omega_{-i}$), there exists a positive constant $D_1$ (respectively, $D_2$) such that $\|f_i(x_i,\mathbf{x}_{-i}) - f_i(y_i,\mathbf{x}_{-i})\|\leq D_1\|x_i-y_i\|$ (respectively, $\|f_i(x_i,\mathbf{x}_{-i}) - f_i(x_i,\mathbf{y}_{-i})\|\leq D_2\|\mathbf{x}_{-i}-\mathbf{y}_{-i}\|$).
\end{Assumption}

\begin{Remark}\label{Remark_existence_NE}
Both Assumptions~\ref{assumption_graph} and \ref{assumption_local_f_lipschitz} are standard and commonly assumed in distributed Nash equilibrium seeking problems. In particular, Assumption~\ref{assumption_local_f_lipschitz} implies that game $\Gamma(N,\{f_i\},\{\Omega_i\},\mathcal{G})$ admits a Nash equilibrium \cite[Prop.2.2]{Debreu1952,Glicksberg1952,Jacquot2018}.
\end{Remark}

\subsection{Preliminaries}
Since the cost function $f_i(x_i,\mathbf{x}_{-i})$ may not be partially differentiable in $x_i$ as in Assumption~\ref{assumption_local_f_lipschitz}, we introduce a Gaussian-smoothed version of the cost function $f_i(x_i,\mathbf{x}_{-i})$ given by \cite{Nesterov2017}
\begin{align*}
  f_{i,\mu^i}(x_i,\mathbf{x}_{-i}) = \frac1{\kappa}\int_{\mathbb{R}} f_i(x_i+\mu^i\xi^i,\mathbf{x}_{-i})e^{-\frac12\|\xi^i\|^2}d\xi^i,
\end{align*}
where $\xi^i \in \mathbb{R}$ is a normally distributed random variable, $\kappa = \int_{\mathbb{R}} e^{-\frac12\|\xi^i\|^2}d\xi^i = (2\pi)^{1/2}$, and $\mu^i \geq 0$ is a smoothing parameter of function $f_{i,\mu^i}(x_i,\mathbf{x}_{-i})$. Then, the randomized gradient-free oracle of $f_i(x_i,\mathbf{x}_{-i})$ can be designed as \cite{Nesterov2017}
\begin{align*}
g^i_{\mu^i}(x_i,\mathbf{x}_{-i}) = \frac{f_i(x_i+\mu^i\xi^i,\mathbf{x}_{-i})-f_i(x_i,\mathbf{x}_{-i})}{\mu^i}\xi^i.
\end{align*}
From the results in \cite{Nesterov2017,Pang2017,Pang2018,Pang2020}, it can be easily shown that the functions $g^i_{\mu^i}(x_i,\mathbf{x}_{-i})$ and $f_{i,\mu^i}(x_i,\mathbf{x}_{-i})$ satisfy some properties which are summarized in the following lemma.
\begin{Lemma}\label{lemma:property_f_mu}
Suppose Assumption~\ref{assumption_local_f_lipschitz} holds. The functions $g^i_{\mu^i}(x_i,\mathbf{x}_{-i})$ and $f_{i,\mu^i}(x_i,\mathbf{x}_{-i})$, $\forall i \in \mathcal{V}$ satisfy the following properties:
\begin{enumerate}
\item The function $f_{i,\mu^i}(x_i,\mathbf{x}_{-i})$ is convex in $x_i$ due to the convexity of $f_i$ in $x_i$. Moreover, $f_{i,\mu^i}(x_i,\mathbf{x}_{-i})$ satisfies
\begin{align*}
f_i(x_i,\mathbf{x}_{-i})\leq f_{i,\mu^i}(x_i,\mathbf{x}_{-i})\leq f_i(x_i,\mathbf{x}_{-i}) + {\mu^i}D_1.
\end{align*}
\item The function $f_{i,\mu^i}(x_i,\mathbf{x}_{-i})$ is partially differentiable in $x_i$ and its partial derivative with respect to $x_i$ satisfies
\begin{align*}
\nabla_{x_i} f_{i,\mu^i}(x_i,\mathbf{x}_{-i}) = \mathbb{E}[g^i_{\mu^i}(x_i,\mathbf{x}_{-i})],
\end{align*}
and is Lipschitz continuous in $x_i$ with a constant $L_1 = \max_{i\in\mathcal{V}}\frac{D_1}{\mu^i}$, and Lipschitz continuous in $\mathbf{x}_{-i}$ with a constant $L_2 = \max_{i\in\mathcal{V}}\frac{D_2}{\mu^i}$, \textit{i.e.,}
\begin{align*}
&\|\nabla_{x_i} f_{i,\mu^i}(x_i,\mathbf{x}_{-i}) - \nabla_{x_i} f_{i,\mu^i}(y_i,\mathbf{x}_{-i})\| \leq L_1\|x_i - y_i\|,\\
&\|\nabla_{x_i} f_{i,\mu^i}(x_i,\mathbf{x}_{-i}) - \nabla_{x_i} f_{i,\mu^i}(x_i,\mathbf{y}_{-i})\| \\
&\quad\quad\quad\quad\quad\quad\quad\quad\quad\quad\quad\quad\quad\quad\quad\leq L_2\|\mathbf{x}_{-i} - \mathbf{y}_{-i}\|.
\end{align*}
Further, $\nabla_{x_i} f_{i,\mu^i}(x_i,\mathbf{x}_{-i})$ always belongs to some $\varepsilon$-subdifferential of function $f_i(x_i,\mathbf{x}_{-i})$, \textit{i.e.,}
\begin{align*}
\nabla_{x_i} f_{i,\mu^i}(x_i,\mathbf{x}_{-i})\in \partial^\varepsilon_{x_i} f_i(x_i,\mathbf{x}_{-i}),\quad \varepsilon = \mu^iD_1.
\end{align*}
Specifically, we have $\nabla_{x_i} f_{i,0}(x_i,\mathbf{x}_{-i}) \in \partial_{x_i} f_i(x_i,\mathbf{x}_{-i})$, when $\mu^i$ tends to 0\footnote{In this paper, we slightly abuse the notation $\mu^i$ to represent the sequence $\mu^i_k$ just for easy presentation without the loss of generality. We mean $\mu^i$ tending to 0 by $\lim_{k\to\infty}\mu^i_k = 0$.}.
\item The random gradient-free oracle $g^i_{\mu^i}(x_i,\mathbf{x}_{-i})$ satisfies
\begin{align*}
\mathbb{E}[\|g^i_{\mu^i}(x_i,\mathbf{x}_{-i})\|] \leq \sqrt{\mathbb{E}[\|g^i_{\mu^i}(x_i,\mathbf{x}_{-i})\|^2]} \leq \mathcal{B},
\end{align*}
where $\mathcal{B} = \sqrt{n+4}D_1$, $n$ is the dimension of $x_i$.
\end{enumerate}
\end{Lemma}

We formulate a smoothed version of game $\Gamma$ with the Gaussian-smoothed cost function $f_{i,\mu^i}$, denoted by $\Gamma_\mu(N,\{f_{i,\mu^i}\},\{\Omega_i\},\mathcal{G})$.
Under Assumption~\ref{assumption_local_f_lipschitz}, the hold of Lemma~\ref{lemma:property_f_mu} implies that the smoothed cost functions $f_{i,\mu^i}$ have similar properties to $f_i$ stated in Assumption~\ref{assumption_local_f_lipschitz}. Hence, game $\Gamma_\mu(N,\{f_{i,\mu^i}\},\{\Omega_i\},\mathcal{G})$ admits a Nash equilibrium for the same reasoning as in Remark~\ref{Remark_existence_NE}.
The following result shows the equivalence of games $\Gamma$ and $\Gamma_\mu$ under certain conditions.
\begin{Lemma}\label{lemma:game_equivalence}
Suppose Assumption~\ref{assumption_local_f_lipschitz} holds. Games $\Gamma$ and $\Gamma_\mu$ are equivalent and share the same Nash equilibria when the smoothing parameter $\mu^i, \forall i \in \mathcal{V}$ tends to 0.
\end{Lemma}
\begin{Proof}
From Remark~\ref{Remark_existence_NE} and previous discussion, Assumption~\ref{assumption_local_f_lipschitz} implies the existence of Nash equilibrium in both games $\Gamma$ and $\Gamma_\mu$. Moreover, applying Squeeze Theorem to Lemma~\ref{lemma:property_f_mu}-1), we have
\begin{align*}
  \lim_{\mu^i\to0}f_{i,\mu^i}(x_i,\mathbf{x}_{-i}) = f_i(x_i,\mathbf{x}_{-i}), \quad \forall i\in \mathcal{V}.
\end{align*}
Then, games $\Gamma$ and $\Gamma_\mu$ share the same number of players, cost functions, action sets and communication graph. Hence the result holds.
\end{Proof}

Next, we make some definitions on the game mappings of games $\Gamma(N,\{f_i\},\{\Omega_i\},\mathcal{G})$ and $\Gamma_\mu(N,\{f_{i,\mu^i}\},\{\Omega_i\},\mathcal{G})$. For game $\Gamma$, since the cost function $f_i$ in game $\Gamma$ is not necessarily differentiable, so the game mapping of game $\Gamma$ refers to a set-valued map $\mathbf{F}(\mathbf{x})$, which is defined as the map of the subdifferentials of all players' cost functions:
\begin{align*}
\mathbf{F}(\mathbf{x})= \prod_{i\in\mathcal{V}}\partial_{x_i}f_i(x_i,\mathbf{x}_{-i}).
\end{align*}
If the cost function $f_i$ is differentiable, then the game mapping $\mathbf{F}(\mathbf{x})$ reduces to a single-valued map. For game $\Gamma_\mu$, since the cost function $f_{i,\mu^i}$ is differentiable, so the game mapping of game $\Gamma_\mu$ refers to a single-valued map $F_\mu$, which is defined by stacking the partial derivatives of all smoothed cost functions:
\begin{align*}
F_\mu(\mathbf{x}) = [\nabla_{x_1}f_{1,\mu^1}(x_1,\mathbf{x}_{-1}),\ldots,\nabla_{x_N}f_{N,\mu^N}(x_N,\mathbf{x}_{-N})]^\top.
\end{align*}
When $\mu^i$ tends to 0 for all $i\in\mathcal{V}$, we denote the game mapping $F_\mu$ by $F_0(\mathbf{x})$, \textit{i.e.},
\begin{align*}
F_0(\mathbf{x}) = \lim_{\mu^i\to0,\forall i\in\mathcal{V}} F_\mu(\mathbf{x})=[\ldots,\nabla_{x_i}f_{i,0}(x_i,\mathbf{x}_{-i}),\ldots]^\top.
\end{align*}
Thus, based on Lemma~\ref{lemma:property_f_mu}-2), we have $F_0(\mathbf{x}) \in \mathbf{F}(\mathbf{x})$.

\section{Gradient-Free Distributed NE Seeking}\label{sec:distr_opt}

In this section, we describe our proposed distributed Nash equilibrium (NE) seeking algorithm in details. 

At time $k$, each player $i\in\mathcal{V}$ maintains an estimate of all players' actions, denoted by $\mathbf{y}^i_k = [y^i_{1,k},\ldots,y^i_{N,k}]^\top \in \mathbb{R}^N$, where $y^i_{j,k}, j\in\mathcal{V}$ represents player $i$'s estimate of player $j$'s action. Hence, at time $k$, every player $l\in\mathcal{V}$ passes its estimate of all players' actions $y^l_{j,k}, j\in\mathcal{V}$ and its own action $x_{l,k}$ to its out-neighbors. Then, for each player $i\in\mathcal{V}$, on receiving the information from its in-neighbors, it updates its own action and the estimate of all players' actions (including the estimate of its own local action) based on the following updating laws:
\begin{subequations}
\begin{align}
x_{i,k+1} &= \mathcal{P}_{\Omega_i}[x_{i,k}-\alpha_kg^i_{\mu^i}(\mathbf{y}^i_k)],\label{eq:update_x}\\
y^i_{j,k+1} &= \sum_{l=1}^N [A]_{il}y^l_{j,k} + \delta_i [A]_{ij}(x_{j,k}-y^i_{j,k}), \quad j\in\mathcal{V}\label{eq:update_y}
\end{align}
\end{subequations}
where $g^i_{\mu^i}(\mathbf{y}^i_k)$ is the randomized gradient-free oracle
\begin{align}
g^i_{\mu^i}(\mathbf{y}^i_k) = \frac{f_i(y^i_{i,k}+\mu^i\xi^i_k,\mathbf{y}^i_{-i,k})-f_i(y^i_{i,k},\mathbf{y}^i_{-i,k})}{\mu^i}\xi^i_k, \label{grad_oracle}
\end{align}
The parameter $\delta_i > 0$ is a constant, and $\alpha_k \geq 0$ is a step-size sequence.
The initial values $x_{i,0}$ and $y^i_{j,0}$ for $i,j\in\mathcal{V}$ can be any real numbers.
The adjacency matrix $A$ is doubly-stochastic as supposed in Assumption~\ref{assumption_graph}.
It should be noted that the design of a doubly-stochastic adjacency matrix $A$ for a given directed graph is non-trivial. The detailed procedures can be referred to the work in \cite{Gharesifard2012}, where two distributed strategies (imbalance-correcting algorithm and load-pushing algorithm) have been developed to construct such matrix under different conditions. With the well-constructed matrix $A$, each player $i$ selects the parameter $\delta_i$ such that $0\leq \delta_i[A]_{ij} < 2[A]_{ii}$ for all $j\in\mathcal{V}$.
The above mentioned procedures are summarized in Algorithm~\ref{algo:rgf_d_dgd}.

\begin{algorithm}
\caption{Gradient-free distributed NE seeking}\label{algo:rgf_d_dgd}
\begin{algorithmic}[1]
\Initialize {arbitrarily generate $x_{i,0},y^i_{j,0},j\in\mathcal{V}$ in $\mathbb{R}$\\ randomly generate $\{\xi^i_k\}_{k\geq0}\sim\mathcal{N}(0,1)$\\ set $\delta_i$ such that $0\leq \delta_i[A]_{ij} < 2[A]_{ii}$, $\forall j\in\mathcal{V}$}
\Iteration {compute $g^i_{\mu^i}(\mathbf{y}^i_k)$ based on (\ref{grad_oracle})\\update variables $x_{i,k+1}$ based on (\ref{eq:update_x})\\update variables $y^i_{j,k+1}$ based on (\ref{eq:update_y})}
\Output {$x_{i,k} \to x^\star$}
\end{algorithmic}
\end{algorithm}
\begin{Remark}
In \eqref{eq:update_y}, it should be noted that $[A]_{ij} = 0$ if player $j$ is NOT an in-neighbor of player $i$, which implies that player $i$ updates the estimate of player $j$'s action only based on the estimates $y^l_{j,k}$ from its in-neighbors $l\in\mathcal{N}^{\text{in}}_i$. On the other hand, if player $j$ is an in-neighbor of player $i$, then $[A]_{ij} \neq 0$ giving rise to an additional error term $\delta_i [A]_{ij}(x_{j,k}-y^i_{j,k})$ in the update of the estimate on player $j$'s action.
\end{Remark}

\section{Convergence Analysis} \label{sec:conv_analysis}
In this section, we study the convergence of the algorithm to the Nash equilibrium for the scenarios of diminishing step-size and constant step-size, respectively. We let $\mathcal{F}_k$ denote the $\sigma$-field generated by the entire history of the random variables from step 0 to $k-1$, \textit{i.e.,}
\begin{align*}
\mathcal{F}_k = \begin{cases}\{x_{i,0}, y^i_{j,0}, i,j \in \mathcal{V}\},&k=0,\\\{x_{i,0}, y^i_{j,0}, \xi^i_s, i,j \in \mathcal{V}; 0\leq s\leq k-1\},&k\geq1.\end{cases}
\end{align*}
Next, we introduce an important property related to the adjacency matrix $A$ summarized in the following lemma:
\begin{Lemma}\label{lemma:A_matrix}
Suppose Assumption~\ref{assumption_graph} holds. Let $\delta_l>0, l \in \mathcal{V}$ be selected such that $0\leq \delta_l[A]_{li} < 2[A]_{ll}$ for all $i\in\mathcal{V}$, where $A$ is the adjacency matrix. Then, there exists a constant $C>0$ and $\gamma \in (0,1)$ such that, for $k\geq 1$, 
the matrix $\tilde{A}_i$ given by
\begin{align*}
[\tilde{A}_i]_{lm} = \begin{cases}[A]_{lm} & \text{if } l\neq m \\|[A]_{ll} -\delta_l[A]_{li}|& \text{if } l= m\end{cases}
\end{align*}
holds that $\|\tilde{A}^k_i\|_\infty<C\gamma^k$.
\end{Lemma}
\begin{Proof}
We first show that all the row sums of $\tilde{A}_i$ are always less than or equal to 1. For any $l\in\mathcal{N}$, if $0\leq \delta_l[A]_{li} \leq [A]_{ll}$, since $A$ is doubly-stochastic, we have $\sum_{n=1}^N[\tilde{A}_i]_{ln} = 1 - \delta_l[A]_{li}\leq 1$; if $[A]_{ll}< \delta_l[A]_{li} < 2[A]_{ll}$, similarly, we have $\sum_{n=1}^N[\tilde{A}_i]_{ln} = 1 + \delta_l[A]_{li} - 2[A]_{ll}< 1$. Hence, we always have $\sum_{n=1}^N[\tilde{A}_i]_{ln} \leq 1$ for any $l\in\mathcal{N}$ where the equal sign holds only if $[A]_{li}=0$.

Next, we show that all the eigenvalues of $\tilde{A}_i$ have magnitude less than or equal to 1.
Let $\lambda$ be an eigenvalue of the matrix $\tilde{A}_i$, and let $\mathbf{v} = [v_1,\ldots,v_N]^\top$ be a corresponding eigenvector. Then we have $\lambda\mathbf{v}=\tilde{A}_i\mathbf{v}$, \textit{i.e.}, for each row $j\in\mathcal{V}$
\begin{align*}
\lambda v_j=[\tilde{A}_i]_{j1}v_1+[\tilde{A}_i]_{j2}v_2+\cdots+[\tilde{A}_i]_{jN}v_N.
\end{align*}
Suppose the $k$-th entry of $\mathbf{v}$ has the maximal absolute value (denoted by $|\hat{v}|$) among all $|v_j|$, $j\in\mathcal{V}$. Then, letting $j = k$ in the above equation, and noting that all the entries of $\tilde{A}_i$ are non-negative and the row sums are less than or equal to 1, we have
\begin{align}
|\lambda||\hat{v}|=|\lambda| |v_k|&=\bigg|\sum_{n=1}^N[\tilde{A}_i]_{kn}v_n\bigg|\nonumber\\
&\leq\bigg(\sum_{n=1}^N[\tilde{A}_i]_{kn}\bigg)|\hat{v}|\leq |\hat{v}|, \label{eq:two_inequality}
\end{align}
which leads to $|\lambda|\leq 1$ as $|\hat{v}|>0$.

Next, we show $|\lambda|\neq 1$ by contradiction. Suppose $|\lambda|= 1$, then the relation \eqref{eq:two_inequality} is true if the equal signs in both inequalities are satisfied, which implies the following properties:
\begin{enumerate}
\item (first equal sign) if $[\tilde{A}_i]_{kn} \neq 0$, then $|v_n|=|\hat{v}|$; 
\item (second equal sign) $[A]_{ki} = 0$.
\end{enumerate}
If player $k$ has an in-neighbour, say player $n\neq k$, then $[\tilde{A}_i]_{kn} = [A]_{kn} \neq 0$. From property 1), we have $|v_n|=|\hat{v}|$. Thus, the $n$-th entry of $\mathbf{v}$ also has the maximal absolute value. That means if player $n\in\mathcal{V}$ is a direct in-neighbor of the player $k\in\mathcal{V}$ with $|v_k|=|\hat{v}|$, then $|v_n|=|\hat{v}|$. Since the graph is strongly connected, thus we can always find a path for each player $n\in\mathcal{V}$ linking to player $k$, \textit{i.e.}, $n\to\cdots\to k$. Thus, from the above analysis, we have $|v_n| = \cdots = |v_k|=|\hat{v}|$ along this path. Therefore, we have $|v_n| =|\hat{v}|$ for all $n\in\mathcal{V}$. From property 2), $|v_n| =|\hat{v}|$ for all $n\in\mathcal{V}$ implies that $[A]_{ni} = 0$ for all $n\in\mathcal{V}$, which is impossible due to the strong connectivity of the graph. Therefore, the eigenvalues of $\tilde{A}_i$ can only have magnitude strictly less than 1, \textit{i.e.}, $|\lambda|<1$.

Finally, we represent $\tilde{A}^k_i$ in the Jordan canonical form for some $P_m$, $J_m$ and $Q_m$. Since all the eigenvalues of $\tilde{A}_i$ have magnitude smaller than 1, then the diagonal entries in $J_m$ are smaller than 1, for all $m$. Thus, there exists a constant $C>0$ and $\gamma \in (0,1)$ such that
 \begin{align*}
\|\tilde{A}^k_i\| &= \bigg\|\sum_{m=1}^NP_m J^k_m Q_m\bigg\|\leq\sum_{m=1}^N\|P_m\|\|Q_m\|\|J^k_m\|\leq C\gamma^k,
\end{align*}
which completes the proof.
\end{Proof}
\begin{Remark}
Similar to the result in \cite[Corollary~1]{Nedic2008}, constants $C$ and $\gamma$ in Lemma~\ref{lemma:A_matrix} depend on the minimum weight (denoted by $\phi$) that each player gives to its own value and the values of its neighbors (\textit{i.e.}, if $[A]_{ij}>0$, then $[A]_{ij}\geq\phi$), the number of players $N$, and the parameters $\{\delta_l\}_{l\in\mathcal{V}}$ selected by all players. Moreover, for larger $N$ and smaller $\phi$, constant $\gamma$ gets closer to 1, implying a slower convergence rate; for $\delta_l \to 0$ or $\delta_l [A]_{li} \to 2[A]_{ll}$, $\forall l\in\mathcal{V}$, matrix $\tilde{A}_i$ reduces to $A$ and constant $\gamma$ gets closer to 1, implying a slower convergence rate.
\end{Remark}

\subsection{Diminishing Step-Size}

In this part, we adopt the diminishing step-size sequence in the proposed algorithm, \textit{i.e.}, the step-size sequence $\alpha_k$ satisfies that $\sum_{k=0}^\infty \alpha_k = \infty$ and $\sum_{k=0}^\infty \alpha^2_k < \infty$. 

Now, we present the result on the consensus property: for any $i\in\mathcal{V}$, each player $l$'s estimate of player $i$'s action $y^l_{i,k}, l\in\mathcal{V}$ converges to player $i$'s real action $x_{i,k}$ as $k$ goes to infinity, which is formally stated in the following theorem.

\begin{Theorem}\label{theorem:consensus}
Suppose Assumptions~\ref{assumption_graph} and \ref{assumption_local_f_lipschitz} hold. Let $\{x_{i,k}\}_{k\geq0}, \{y^l_{i,k}\}_{k\geq0}$, $i,l\in\mathcal{V}$ be the sequences generated by \eqref{eq:update_x} and \eqref{eq:update_y}, respectively, with a step-size sequence $\{\alpha_k\}_{k\geq0}$ satisfying $\sum_{k=0}^\infty \alpha_k = \infty$ and $\sum_{k=0}^\infty \alpha^2_k < \infty$, and a positive constant $\delta_l, l\in\mathcal{V}$ satisfying $0\leq \delta_l[A]_{li} < 2[A]_{ll}$ for all $i\in\mathcal{V}$, where $A$ is the adjacency matrix. Then, we have
 \begin{align*}
\lim_{k\to\infty}\mathbb{E}[\|x_{i,k}-y^l_{i,k}\|] = 0,\quad i,l\in\mathcal{V}.
\end{align*}
\end{Theorem}
\begin{Proof}
It can be obtained from \eqref{eq:update_x} and \eqref{eq:update_y} that
\begin{align*}
x_{i,k+1} &= \mathcal{P}_{\Omega_i}[x_{i,k}-\alpha_{k}g^i_{\mu^i}(\mathbf{y}^i_k)],\\
y^l_{i,k+1} &= \sum_{m=1}^N[A]_{lm}y^m_{i,k}+\delta_{l}[A]_{li}(x_{i,k}-y^l_{i,k}).
\end{align*}
Then, taking the subtraction and applying the norm
\begin{align*}
&\|x_{i,k+1}-y^l_{i,k+1}\|\\ 
\leq &\bigg\|\sum_{m=1}^N[A]_{lm}(x_{i,k}-y^m_{i,k})-\delta_l[A]_{li}(x_{i,k}-y^l_{i,k})\bigg\|\\
&+\alpha_{k}\|g^i_{\mu^i}(\mathbf{y}^i_k)\|\\
=&\bigg\|\sum_{m=1}^N[\tilde{A}_i]_{lm}(x_{i,k}-y^m_{i,k})\bigg\|+\alpha_{k}\|g^i_{\mu^i}(\mathbf{y}^i_k)\|\\
\leq &\sum_{m=1}^N[\tilde{A}_i]_{lm}\|x_{i,k}-y^m_{i,k}\|+\alpha_{k}\|g^i_{\mu^i}(\mathbf{y}^i_k)\|,
\end{align*} 
where the first inequality follows from the projection's non-expansive property, and the equality holds by the defintion of $\tilde{A}_i$ as in Lemma~\ref{lemma:A_matrix}. Hence, we obtain that
\begin{align*}
\|x_{i,k}-y^l_{i,k}\|\leq& \sum_{m=1}^N[\tilde{A}^k_i]_{lm}\|x_{i,0}-y^m_{i,0}\|\\
&+\sum_{r=1}^{k-1}\sum_{m=1}^N[\tilde{A}^{k-r}_i]_{lm}\alpha_{r-1}\|g^i_{\mu^i}(\mathbf{y}^i_{r-1})\|\\
&+\alpha_{k-1}\|g^i_{\mu^i}(\mathbf{y}^i_{k-1})\|.
\end{align*}
Then, taking the total expectation, it follows from Lemmas~\ref{lemma:property_f_mu}-3) and \ref{lemma:A_matrix} that
\begin{align}
\mathbb{E}[\|x_{i,k}-y^l_{i,k}\|]&\leq NC\hat{\sigma}\gamma^k \nonumber\\
&\quad+NC\mathcal{B}\sum_{r=1}^{k-1}\gamma^{k-r}\alpha_{r-1}+\mathcal{B}\alpha_{k-1}, \label{eq:x_i_minus_y^l_i}
\end{align}
where $\hat{\sigma} = \max_{i\in\mathcal{V}}\sigma_i$ and $\sigma_i = \max_{m\in\mathcal{V}}|x_{i,0}-y^m_{i,0}|$.
Taking the limit $k\to\infty$ and noting that $\lim_{k\to\infty}\alpha_k=0$, the desired result follows from \cite[Lemma 4-1)]{Pang2017}.
\end{Proof}
\begin{Remark}
Theorem~\ref{theorem:consensus} is a characterization of the consensus property of the algorithm. For any $i\in\mathcal{V}$, each player $l$'s estimate of player $i$'s action $y^l_{i,k}, l\in\mathcal{V}$ converges to player $i$'s real action $x_{i,k}$ as $k$ goes to infinity.
\end{Remark}

Now, we make a mild assumption on the uniqueness of the Nash equilibrium in game $\Gamma$ as follows.
\begin{Assumption}\label{assumption_game_mapping}
The game mapping $\mathbf{F}$ of game $\Gamma$ is strictly monotone on $\Omega$, \textit{i.e.}, for any $\mathbf{x},\mathbf{y}\in\Omega$, $\mathbf{x}\neq\mathbf{y}$, $\mathbf{f}(\mathbf{x}), \mathbf{f}(\mathbf{y}) \in \mathbf{F}$, we have $\langle \mathbf{f}(\mathbf{x})-\mathbf{f}(\mathbf{y}), \mathbf{x}-\mathbf{y} \rangle> 0$.
\end{Assumption}
\begin{Remark}
Assumption~\ref{assumption_game_mapping} ensures the uniqueness of the Nash equilibrium in game $\Gamma$.
\end{Remark}

Now, we are ready to establish the convergence of all players' actions to the unique Nash equilibrium of game $\Gamma$, which is formally stated in the following theorem.
\begin{Theorem}\label{theorem:optimality}
Suppose Assumptions~\ref{assumption_graph}, \ref{assumption_local_f_lipschitz} and \ref{assumption_game_mapping} hold. Let $\mathbf{x}^\star=(x^\star_i,\mathbf{x}^\star_{-i})$ be the action profile at the unique Nash equilibrium of the game $\Gamma$. Let $\{x_{i,k}\}_{k\geq0}, \{y^l_{i,k}\}_{k\geq0}$, $i,l\in\mathcal{V}$ be the sequences generated by \eqref{eq:update_x} and \eqref{eq:update_y}, respectively, with a step-size sequence $\{\alpha_k\}_{k\geq0}$ satisfying $\sum_{k=0}^\infty \alpha_k = \infty$ and $\sum_{k=0}^\infty \alpha^2_k < \infty$, and a positive constant $\delta_l, l\in\mathcal{V}$ satisfying $0\leq \delta_l[A]_{li} < 2[A]_{ll}$ for all $i\in\mathcal{V}$, where $A$ is the adjacency matrix. Then, the sequence $\{\mathbf{x}_k\}_{k\geq0}$ converges to $\mathbf{x}^\star$ almost surely when the smoothing parameter $\mu^i, \forall i \in \mathcal{V}$ tends to 0.
\end{Theorem}
\begin{Proof}
Noting that there might be multiple Nash equilibria in game $\Gamma_\mu$, we let $\mathbf{x}_\mu^\star=(x^\star_{i,\mu^i},\mathbf{x}^\star_{-i,\mu^{-i}})$ be an action profile at one of them. Applying \cite[Lemma~1]{Salehisadaghiani2016} yields
\begin{align*}
x^\star_{i,\mu^i} = \mathcal{P}_{\Omega_i}[x^\star_{i,\mu^i} - \alpha_{k} \nabla_{x_i} f_{i,\mu^i}(\mathbf{x}^\star_\mu)], \quad i\in \mathcal{V}.
\end{align*}
Thus, subtracting \eqref{eq:update_x} by the above equation, and taking the norm
\begin{align*}
&\|x_{i,k+1}-x^\star_{i,\mu^i}\|^2 \\
\leq &\|(x_{i,k}-x^\star_{i,\mu^i})-\alpha_{k}(g^i_{\mu^i}(\mathbf{y}^i_k)-\nabla_{x_i} f_{i,\mu^i}(\mathbf{x}^\star_\mu))\|^2\\
= &\|x_{i,k}-x^\star_{i,\mu^i}\|^2 + \alpha_{k}^2\|g^i_{\mu^i}(\mathbf{y}^i_k)-\nabla_{x_i} f_{i,\mu^i}(\mathbf{x}^\star_\mu)\|^2\\
&-2\alpha_{k}\langle x_{i,k}-x^\star_{i,\mu^i},g^i_{\mu^i}(\mathbf{y}^i_k)-\nabla_{x_i} f_{i,\mu^i}(\mathbf{x}^\star_\mu)\rangle,
\end{align*}
where we have applied the projection's non-expansive property. Taking the conditional expectation on $\mathcal{F}_{k}$, we obtain
\begin{align}
&\mathbb{E}[\|x_{i,k+1}-x^\star_{i,\mu^i}\|^2|\mathcal{F}_{k}] \nonumber\\
\leq &\|x_{i,k}-x^\star_{i,\mu^i}\|^2 + \alpha_{k}^2\mathbb{E}[\|g^i_{\mu^i}(\mathbf{y}^i_k)-\nabla_{x_i} f_{i,\mu^i}(\mathbf{x}^\star_\mu)\|^2|\mathcal{F}_{k}]\nonumber\\
&-2\alpha_{k}\langle x_{i,k}-x^\star_{i,\mu^i},\nabla_{x_i} f_{i,\mu^i}(\mathbf{y}^i_k)-\nabla_{x_i} f_{i,\mu^i}(\mathbf{x}^\star_\mu)\rangle\nonumber\\
=& \|x_{i,k}-x^\star_{i,\mu^i}\|^2 + \alpha_{k}^2\mathbb{E}[\|g^i_{\mu^i}(\mathbf{y}^i_k)-\nabla_{x_i} f_{i,\mu^i}(\mathbf{x}^\star_\mu)\|^2|\mathcal{F}_{k}]\nonumber\\
&-2\alpha_{k}\langle x_{i,k}-x^\star_{i,\mu^i},\nabla_{x_i} f_{i,\mu^i}(\mathbf{y}^i_k)-\nabla_{x_i} f_{i,\mu^i}(\mathbf{x}_k)\rangle\nonumber\\
&-2\alpha_{k}\langle x_{i,k}-x^\star_{i,\mu^i},\nabla_{x_i} f_{i,\mu^i}(\mathbf{x}_k)-\nabla_{x_i} f_{i,\mu^i}(\mathbf{x}^\star_\mu)\rangle. \label{eq:x_i_minus_x^star_mu}
\end{align}
It is noted that
\begin{align*}
&\mathbb{E}[\|g^i_{\mu^i}(\mathbf{y}^i_k)-\nabla_{x_i} f_{i,\mu^i}(\mathbf{x}^\star_\mu)\|^2|\mathcal{F}_{k}]\\
&\quad\leq 2\mathbb{E}[\|g_{\mu^i}(\mathbf{y}^i_k)\|^2|\mathcal{F}_{k}]+2\|\nabla_{x_i} f_{i,\mu^i}(\mathbf{x}^\star_\mu)\|^2\leq 4\mathcal{B}^2
\end{align*}
and
\begin{align*}
&\quad-2\alpha_{k}\langle x_{i,k}-x^\star_{i,\mu^i},\nabla_{x_i} f_{i,\mu^i}(\mathbf{y}^i_k)-\nabla_{x_i} f_{i,\mu^i}(\mathbf{x}_k)\rangle\\
&=-2\alpha_{k}\langle x_{i,k}-x^\star_{i,\mu^i},\nabla_{x_i} f_{i,\mu^i}(\mathbf{y}^i_k)-\nabla_{x_i} f_{i,\mu^i}(y^i_k,\mathbf{x}_{-i,k})\rangle\\
&\quad-2\alpha_{k}\langle x_{i,k}-x^\star_{i,\mu^i},\nabla_{x_i} f_{i,\mu^i}(y^i_k,\mathbf{x}_{-i,k}) - \nabla_{x_i} f_{i,\mu^i}(\mathbf{x}_{k})\rangle\\
&\leq2L_1\alpha_{k}\|x_{i,k}-x^\star_{i,\mu^i}\|\|y^i_k-x_{i,k}\|\\
&\quad+2L_2\alpha_{k}\|x_{i,k}-x^\star_{i,\mu^i}\|\|\mathbf{y}^i_{-i,k}-\mathbf{x}_{-i,k}\|\\
&\leq L_1(1+\|x_{i,k}-x^\star_{i,\mu^i}\|^2)\alpha_{k}\|y^i_k-x^i_k\|\\
&\quad+L_2(1+\|x_{i,k}-x^\star_{i,\mu^i}\|^2)\alpha_{k}\|\mathbf{y}^i_{-i,k}-\mathbf{x}_{-i,k}\|\\
&\leq\hat{L}(1+\|x_{i,k}-x^\star_{i,\mu^i}\|^2)\bigg(NC\hat{\sigma}\alpha_{k}\gamma^k\\
&\quad+NC\mathcal{B}\sum_{r=1}^{k-1}\gamma^{k-r}\alpha_{k}\alpha_{r-1}+\mathcal{B}\alpha_{k}\alpha_{k-1}\bigg),
\end{align*}
where $\hat{L} = L_1 + \sqrt{N-1}L_2$, we have applied $2\|\mathbf{a}\|\leq 1 + \|\mathbf{a}^2\|$ in the second inequality and \eqref{eq:x_i_minus_y^l_i} in the third inequality. Thus, combining the above results to \eqref{eq:x_i_minus_x^star_mu} and summing over $i\in\mathcal{V}$
\begin{align}
&\mathbb{E}[\|\mathbf{x}_{k+1}-\mathbf{x}^\star_{\mu}\|^2|\mathcal{F}_{k}] \leq  (1+\eta_k)\|\mathbf{x}_{k}-\mathbf{x}^\star_{\mu}\|^2 + w_k \nonumber\\
&\quad\quad\quad\quad\quad\quad\quad- 2\alpha_{k}\langle \mathbf{x}_{k}-\mathbf{x}^\star_{\mu},F_{\mu}(\mathbf{x}_{k})- F_{\mu}(\mathbf{x}^\star_\mu)\rangle, \label{eq:big_x_k_plus_minus_big_x_star_mu}
\end{align}
where
\begin{align*}
\eta_k &= \hat{L}\bigg(NC\hat{\sigma}\alpha_{k}\gamma^k\\
&\quad+NC\mathcal{B}\sum_{r=1}^{k-1}\gamma^{k-r}\alpha_{k}\alpha_{r-1}+\mathcal{B}\alpha_{k}\alpha_{k-1}\bigg)\\
w_k &= 4N\mathcal{B}^2\alpha_{k}^2+\hat{L}\bigg(N^2C\hat{\sigma}\alpha_{k}\gamma^k\\
&\quad+N^2C\mathcal{B}\sum_{r=1}^{k-1}\gamma^{k-r}\alpha_{k}\alpha_{r-1}+N\mathcal{B}\alpha_{k}\alpha_{k-1}\bigg).
\end{align*}
Taking the limit $\mu^i \to \infty$ for $\forall i \in \mathcal{V}$, it follows from Lemma~\ref{lemma:game_equivalence} that $\lim_{\mu^i \to \infty, \forall i \in \mathcal{V}} \mathbf{x}^\star_\mu = \mathbf{x}^\star$. Thus, it can be obtained from \eqref{eq:big_x_k_plus_minus_big_x_star_mu} that
\begin{align}
&\mathbb{E}[\|\mathbf{x}_{k+1}-\mathbf{x}^\star\|^2|\mathcal{F}_{k}] \leq  (1+\eta_k)\|\mathbf{x}_{k}-\mathbf{x}^\star\|^2 + w_k \nonumber\\
&\quad\quad\quad\quad\quad\quad\quad- 2\alpha_{k}\langle \mathbf{x}_{k}-\mathbf{x}^\star,F_0(\mathbf{x}_{k})- F_0(\mathbf{x}^\star)\rangle. \label{eq:big_x_k_plus_minus_big_x_star}
\end{align}
Following the results in \cite[Lemma~3]{Pang2017} and the step-size $\sum_{k=0}^\infty \alpha_k = \infty$, $\sum_{k=0}^\infty \alpha^2_k < \infty$, we have
\begin{align*}
\sum_{k=0}^\infty\eta_k<\infty,\quad
\sum_{k=0}^\infty w_k <\infty.
\end{align*}
Applying Lemma~11 in \cite[Ch.~2]{Polyak1987} to \eqref{eq:big_x_k_plus_minus_big_x_star}, we can obtain that $\|\mathbf{x}_{k}-\mathbf{x}^\star\|$ converges almost surely, and $\sum_{k=0}^\infty\alpha_{k}\langle \mathbf{x}_{k}-\mathbf{x}^\star,F_0(\mathbf{x}_{k})- F_0(\mathbf{x}^\star)\rangle <\infty$. From Assumption~\ref{assumption_game_mapping} and the fact that $F_0(\mathbf{x}) \in \mathbf{F}(\mathbf{x})$, we have $\langle \mathbf{x}_{k}-\mathbf{x}^\star,F_0(\mathbf{x}_{k})- F_0(\mathbf{x}^\star)\rangle \geq 0$. Together with the step-size $\sum_{k=0}^\infty \alpha_k = \infty$, $\sum_{k=0}^\infty \alpha^2_k < \infty$, we obtain
\begin{align*}
\liminf_{k\to\infty} \mathbf{x}_{k}=\mathbf{x}^\star.
\end{align*}
Since $\|\mathbf{x}_{k}-\mathbf{x}^\star\|$ converges almost surely, we obtain the desired result.
\end{Proof}

\begin{Remark}
Theorem~\ref{theorem:optimality} shows that the players' action profile $\mathbf{x}_{k}$ will converge to the unique Nash equilibrium $\mathbf{x}^\star$ of game $\Gamma$ by selecting the diminishing smoothing parameter sequence.
\end{Remark}

\subsection{Constant Step-Size}

In this part, we suppose the step-size $\alpha_k = \alpha$, which is a positive constant.

A similar result to Theorem~\ref{theorem:consensus} on the consensus property can be established. Instead of achieving the exact convergence, for any $i\in\mathcal{V}$, each player $l$'s estimate of player $i$'s action $y^l_{i,k}, l\in\mathcal{V}$ approximately converges to player $i$'s real action $x_{i,k}$ with an error proportional to the step-size. The following theorem formally states the result.

\begin{Theorem}\label{theorem:consensus_2}
Suppose Assumptions~\ref{assumption_graph} and \ref{assumption_local_f_lipschitz} hold. Let $\{x_{i,k}\}_{k\geq0}, \{y^l_{i,k}\}_{k\geq0}$, $i,l\in\mathcal{V}$ be the sequences generated by \eqref{eq:update_x} and \eqref{eq:update_y}, respectively, with a constant step-size sequence $\alpha_k = \alpha$, and a positive constant $\delta_l, l\in\mathcal{V}$ satisfying $0\leq \delta_l[A]_{li} < 2[A]_{ll}$ for all $i\in\mathcal{V}$, where $A$ is the adjacency matrix. Then, we have
 \begin{align*}
\limsup_{k\to\infty}\mathbb{E}[\|x_{i,k}-y^l_{i,k}\|] \leq \bigg(\frac{\gamma NC\mathcal{B}}{1-\gamma}+\mathcal{B}\bigg)\alpha.
\end{align*}
\end{Theorem}
\begin{Proof}
Following same arguments as in Theorem~\ref{theorem:consensus}, the result holds by taking the limsup on both sides of \eqref{eq:x_i_minus_y^l_i}.
\end{Proof}

Next, we introduce a slightly stronger assumption compared to Assumption~\ref{assumption_game_mapping} on the uniqueness of the Nash equilibrium in game $\Gamma$.
\begin{Assumption}\label{assumption_game_mapping_2}
The game mapping $\mathbf{F}$ of game $\Gamma$ is strongly monotone on $\Omega$ with a constant $\chi>0$, \textit{i.e.}, for any $\mathbf{x},\mathbf{y}\in\Omega$, $\mathbf{f}(\mathbf{x}), \mathbf{f}(\mathbf{y}) \in \mathbf{F}$, we have $\langle \mathbf{f}(\mathbf{x})-\mathbf{f}(\mathbf{y}), \mathbf{x}-\mathbf{y} \rangle\geq \chi\|\mathbf{x}-\mathbf{y}\|^2$.
\end{Assumption}
\begin{Remark}
Assumption~\ref{assumption_game_mapping_2} also ensures the uniqueness of the Nash equilibrium in game $\Gamma$.
\end{Remark}

Now, we are ready to characterize the approximate convergence of all players' actions to the Nash equilibrium of game $\Gamma$, which is formally stated in the following theorem.

\begin{Theorem}\label{theorem:optimality_2}
Suppose Assumptions~\ref{assumption_graph}, \ref{assumption_local_f_lipschitz} and \ref{assumption_game_mapping_2} hold. Let $\mathbf{x}^\star=(x^\star_i,\mathbf{x}^\star_{-i})$ be the action profile at the unique Nash equilibrium of game $\Gamma$. Let $\{x_{i,k}\}_{k\geq0}, \{y^l_{i,k}\}_{k\geq0}$, $i,l\in\mathcal{V}$ be the sequences generated by \eqref{eq:update_x} and \eqref{eq:update_y}, respectively, with a constant step-size sequence $\alpha_k=\alpha$ satisfying the following condition
\begin{align}
0<2\chi\alpha - \hat{L}\bigg(\frac{\gamma NC\mathcal{B}}{1-\gamma}+\mathcal{B}\bigg)\alpha^2<1, \label{eq:const_step_condition}
\end{align}
and a positive constant $\delta_l, l\in\mathcal{V}$ satisfying $0\leq \delta_l[A]_{li} < 2[A]_{ll}$ for all $i\in\mathcal{V}$, where $A$ is the adjacency matrix.
Then, with the smoothing parameter $\mu^i, \forall i \in \mathcal{V}$ tending to 0, the sequence $\{\mathbf{x}_k\}_{k\geq0}$ satisfies
\begin{align*}
\limsup_{k\to\infty}\mathbb{E}[\|\mathbf{x}_{k}-\mathbf{x}^\star\|^2]\leq \frac{[4N\mathcal{B}^2+\hat{L}(\frac{\gamma N^2C\mathcal{B}}{1-\gamma}+N\mathcal{B})]\alpha}{2\chi-\hat{L}(\frac{\gamma NC\mathcal{B}}{1-\gamma}+\mathcal{B})\alpha}.
\end{align*}
\end{Theorem}
\begin{Proof}
From Assumption~\ref{assumption_game_mapping_2} and the fact that $F_0(\mathbf{x}) \in \mathbf{F}(\mathbf{x})$, we have $\langle \mathbf{x}_{k}-\mathbf{x}^\star,F_0(\mathbf{x}_{k})- F_0(\mathbf{x}^\star)\rangle \geq \chi\|\mathbf{x}_{k}-\mathbf{x}^\star\|$.
Following the same arguments as in Theorem~\ref{theorem:optimality}, and applying the above results to \eqref{eq:big_x_k_plus_minus_big_x_star}, we obtain
\begin{align*}
\mathbb{E}[\|\mathbf{x}_{k+1}-\mathbf{x}^\star\|^2|\mathcal{F}_{k}] \leq  (1+\eta_k- 2\chi\alpha)\|\mathbf{x}_{k}-\mathbf{x}^\star\|^2 + w_k,
\end{align*}
where
\begin{align*}
\eta_k &= \hat{L}\bigg(NC\hat{\sigma}\alpha\gamma^k+NC\mathcal{B}\alpha^2\sum_{r=1}^{k-1}\gamma^{k-r}+\mathcal{B}\alpha^2\bigg)\\
w_k &= 4N\mathcal{B}^2\alpha^2+\hat{L}\bigg(N^2C\hat{\sigma}\alpha\gamma^k\\
&\quad\quad\quad\quad\quad+N^2C\mathcal{B}\alpha^2\sum_{r=1}^{k-1}\gamma^{k-r}+N\mathcal{B}\alpha^2\bigg).
\end{align*}
Taking the total expectation and the limsup on both sides, we complete the proof based on the step-size condition \eqref{eq:const_step_condition}.
\end{Proof}

\begin{Remark}
In general, if the constant step-size $\alpha$ is set small, then the step-size condition \eqref{eq:const_step_condition} can be satisfied.
Theorem~\ref{theorem:optimality_2} shows that all players' actions approximately converge to the Nash equilibrium of game $\Gamma$ with an error depending on the step-size $\alpha$, the number of players $N$, the cost function parameters $D_1, D_2$ ($\hat{L}$ and $\mathcal{B}$ are functions of $D_1, D_2$) and the communication topology $\gamma$. It should also be noted that if the step-size $\alpha$ is small, then the error bound is close to 0.
\end{Remark}

\section{Numerical Simulations}\label{sec:simulation}

In this section, we demonstrate the performance of the proposed algorithm by a numerical example. Consider an energy consumption game of $N$ players for Heating Ventilation and Air Conditioning (HVAC) system (see \cite{Ye2017}), where the cost function of each player $i$ can be modeled by the following quadratic function:
\begin{equation*}
f_i(\mathbf{x}) = a_i(x_i - x^r_i)^2 + \bigg(b\sum_{j=1}^N x_j + c\bigg)x_i, \quad x_i\in\Omega_i,
\end{equation*}
where $a_i>0, b>0, c$ and $x^r_i$ are constants for $i\in\mathcal{V}$. It is easy to verify that Assumptions~2, 3 and 4 are satisfied. Throughout the simulation, we let $a_i = 1$ for $i\in\mathcal{V}$, $b = 0.1$ and $c = 10$. In the following simulation, we investigate the effectiveness of the proposed algorithm from the perspectives of network topology and number of players, followed by a comparison with the gradient-based counterpart.

\subsection{Network Topology}
\begin{figure}[!t]
    \centering
    \subfloat[$\mathcal{G}_1$]{{\includegraphics[width=1.1in]{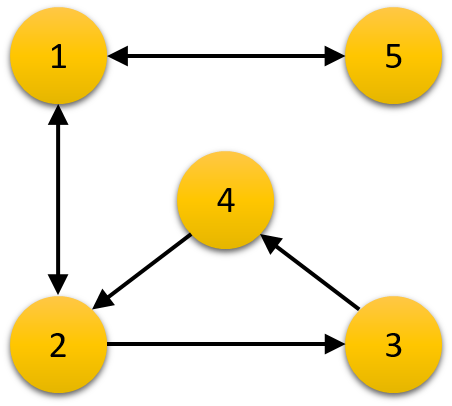} }}%
    \subfloat[$\mathcal{G}_2$]{{\includegraphics[width=1.1in]{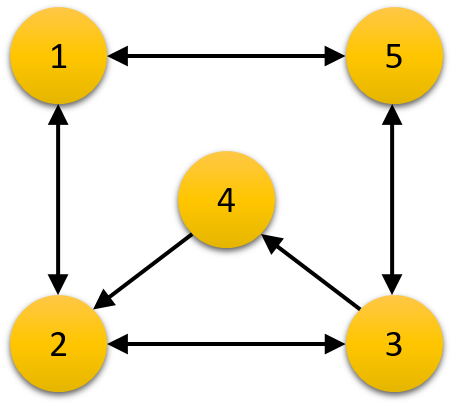} }}%
    \subfloat[$\mathcal{G}_3$]{{\includegraphics[width=1.1in]{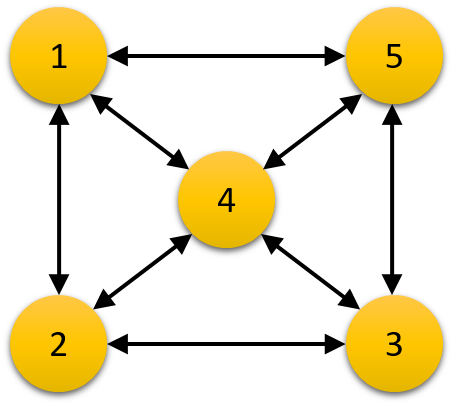} }}%
    \caption{Three different network topologies.}%
    \label{fig:network.png}%
\end{figure}

\begin{figure}[!t]
\centering
\includegraphics[width=3.4in]{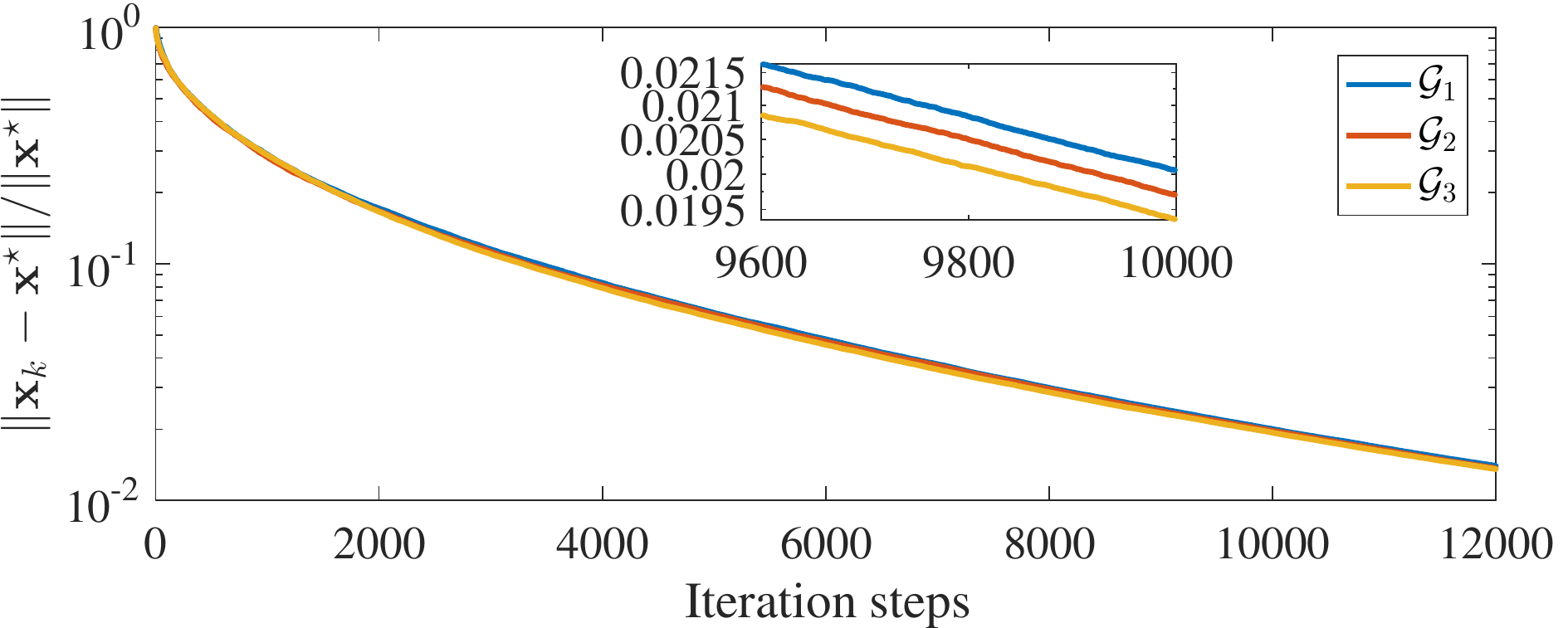}  %
\caption{The relative errors of all players' actions produced by the proposed gradient-free method ($\alpha_k= 0.1/\sqrt{k+1}$) under three different communication networks.}
\label{fig:network_compare.PNG}
\end{figure}

In this part, we first consider $N=5$ players under three different communication graphs as shown in Fig.~\ref{fig:network.png}. Obviously, all these digraphs are strongly connected, hence Assumption~1 is satisfied.
Constant $x^r_i$ for $i\in\{1,\ldots,5\}$ is set to $10, 15, 20, 25, 30$, respectively. 
For the implementation of the algorithm, the step-size $\alpha_k$ is set to $0.1/\sqrt{k+1}$. Besides, we let the smoothing parameter sequence to be diminishing, \textit{e.g.}, $\mu^i_k = 10^{-2}/{k+1}$ and $\delta_i = 0.5$ for $i\in\{1,\ldots,5\}$. The initial values of all players' actions $\mathbf{x}_0$ and the estimates of all players' actions $\mathbf{y}^i_0$ from player $i$ for $i\in\{1,\ldots,5\}$ are all set to 0. The relative errors of all players' actions ($\|\mathbf{x}_k-\mathbf{x}^\star\|/\|\mathbf{x}^\star\|$) produced by the proposed gradient-free method with diminishing step-size for three different network topologies are plotted in Fig.~\ref{fig:network_compare.PNG}.
As can be observed, convergence can be achieved under all three network topologies. Specifically, the performance is better for the graph with more edges due to the increased number of communication channels.

\subsection{Number of Players}
\begin{figure}[!t]
    \centering
    \includegraphics[width=2.3in]{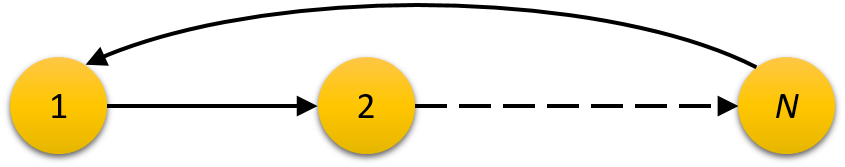}%
    \caption{Network topology for $N$ players.}%
    \label{fig:network_number_of_players.PNG}%
\end{figure}

\begin{figure}[!t]
\centering
\includegraphics[width=3.4in]{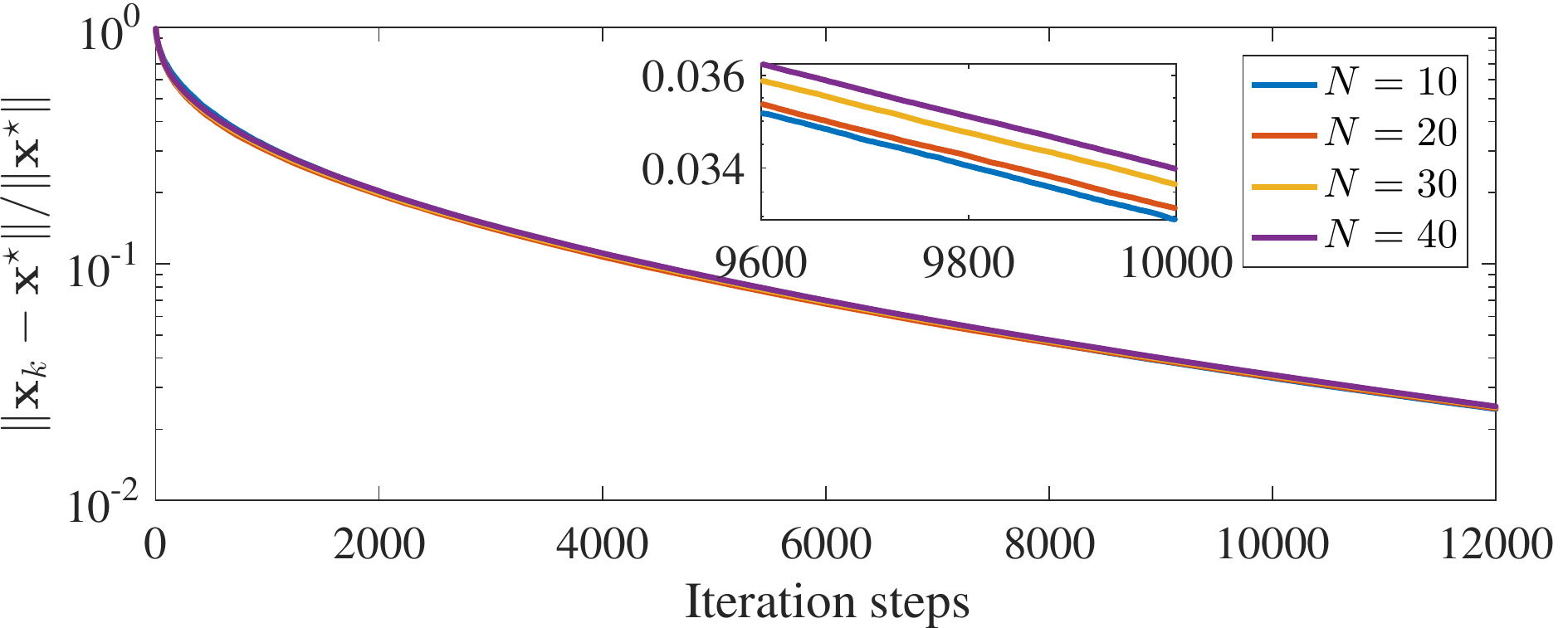}  %
\caption{The relative errors of all players' actions produced by the proposed gradient-free method ($\alpha_k= 0.1/\sqrt{k+1}$) with different number of players.}
\label{fig:numAgent_compare.PNG}
\end{figure}

In this part, we increase the number of players to $N = 10, 20, 30$ and $40$ under a strongly connected communication graph as shown in Fig.~\ref{fig:network_number_of_players.PNG}. 
We set $x^r_i = 2i$ for $i\in\{1,\ldots,N\}$. The rest of parameters are set the same as in section~\ref{sec:simulation}-A. It is shown in Fig.~\ref{fig:numAgent_compare.PNG} that the relative errors of all players' actions ($\|\mathbf{x}_k-\mathbf{x}^\star\|/\|\mathbf{x}^\star\|$) produced by the proposed gradient-free method with diminishing step-size for $N = 10, 20, 30$ and $40$.
As can be seen, the algorithm is scalable to different number of players, and the convergence result is better for smaller number of players, which is as expected. 

\subsection{Gradient-Free vs. Gradient-Based Algorithm}
\begin{figure}[!t]
\centering
\includegraphics[width=3.4in]{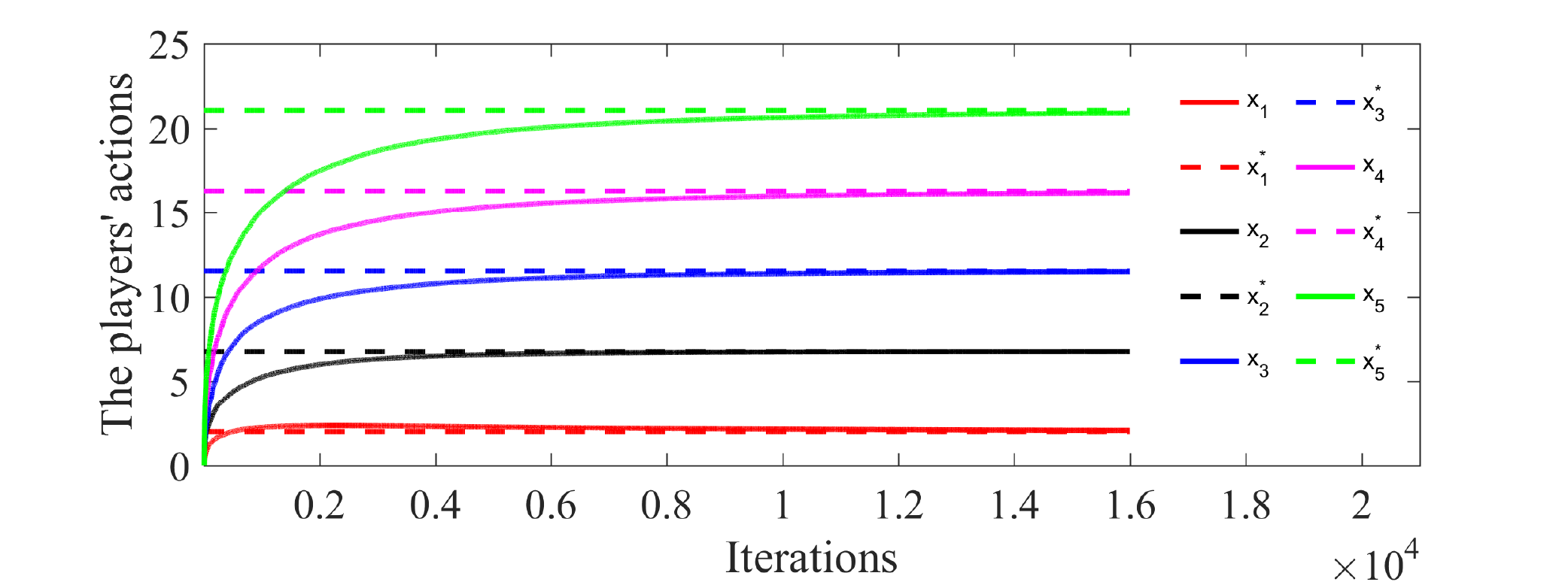}  %
\caption{The plot of all players' actions produced by the proposed gradient-free method with diminishing step-size $\alpha_k= 0.1/\sqrt{k+1}$.}
\label{fig:diminishing_stepsize_rgf.PNG}
\end{figure}

\begin{figure}[!t]
\centering
\includegraphics[width=3.4in]{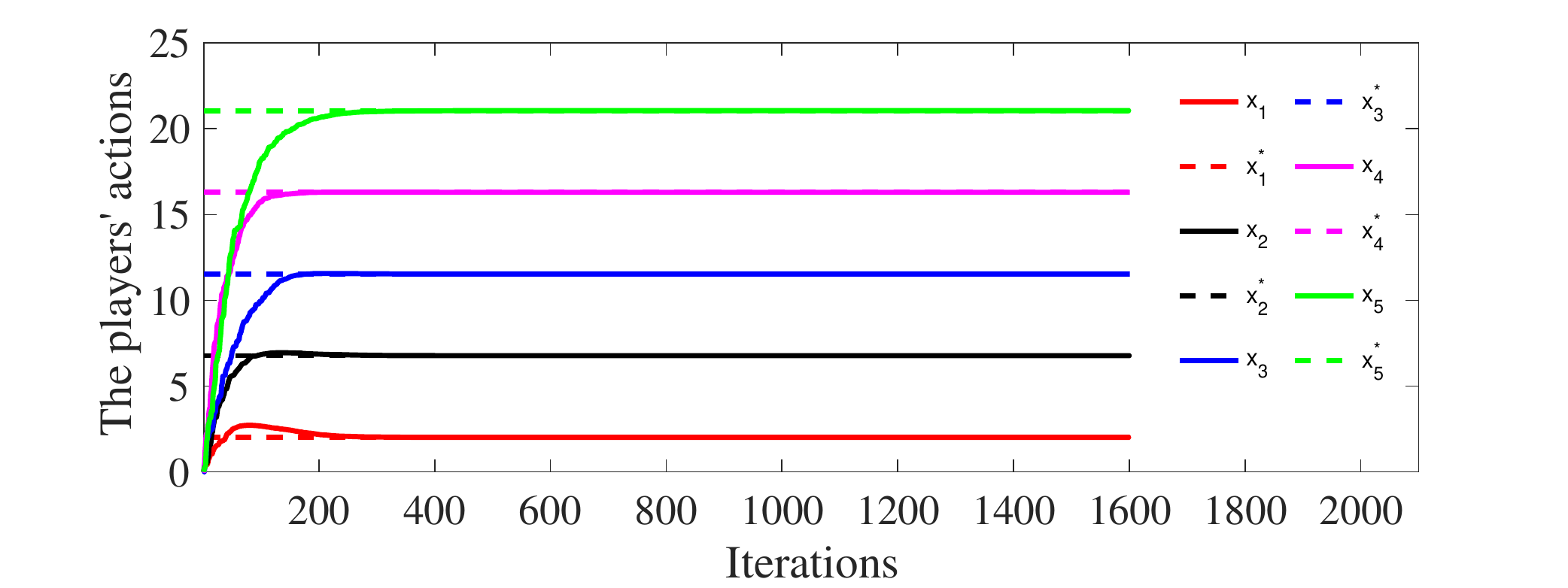} %
\caption{The plot of all players' actions produced by the proposed gradient-free method with constant step-size $\alpha_k= 0.1$.}
\label{fig:constant_stepsize_rgf.PNG}
\end{figure}

\begin{figure}[!t]
\centering
\includegraphics[width=3.4in]{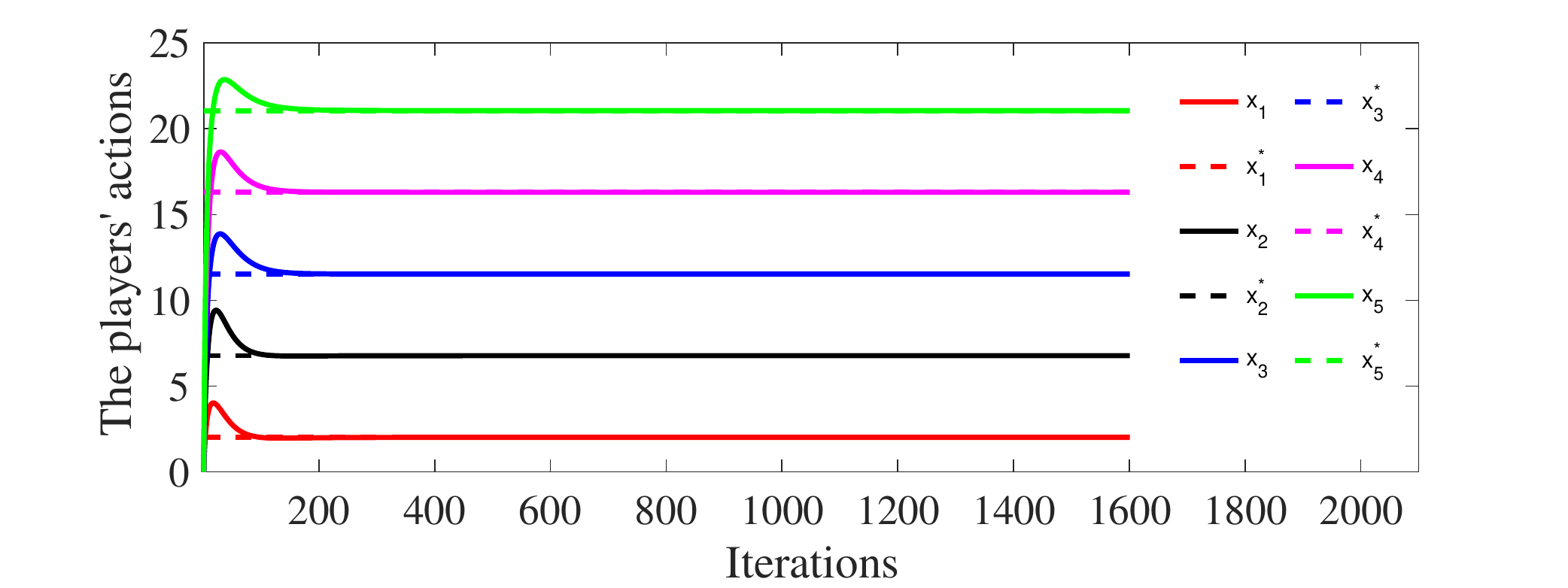}  %
\caption{The plot of all players' actions produced by the gradient-based method with diminishing step-size $\alpha_k= 0.1/\sqrt{k+1}$.}
\label{fig:diminishing_stepsize_grad.PNG}
\end{figure}

\begin{figure}[!t]
\centering
\includegraphics[width=3.4in]{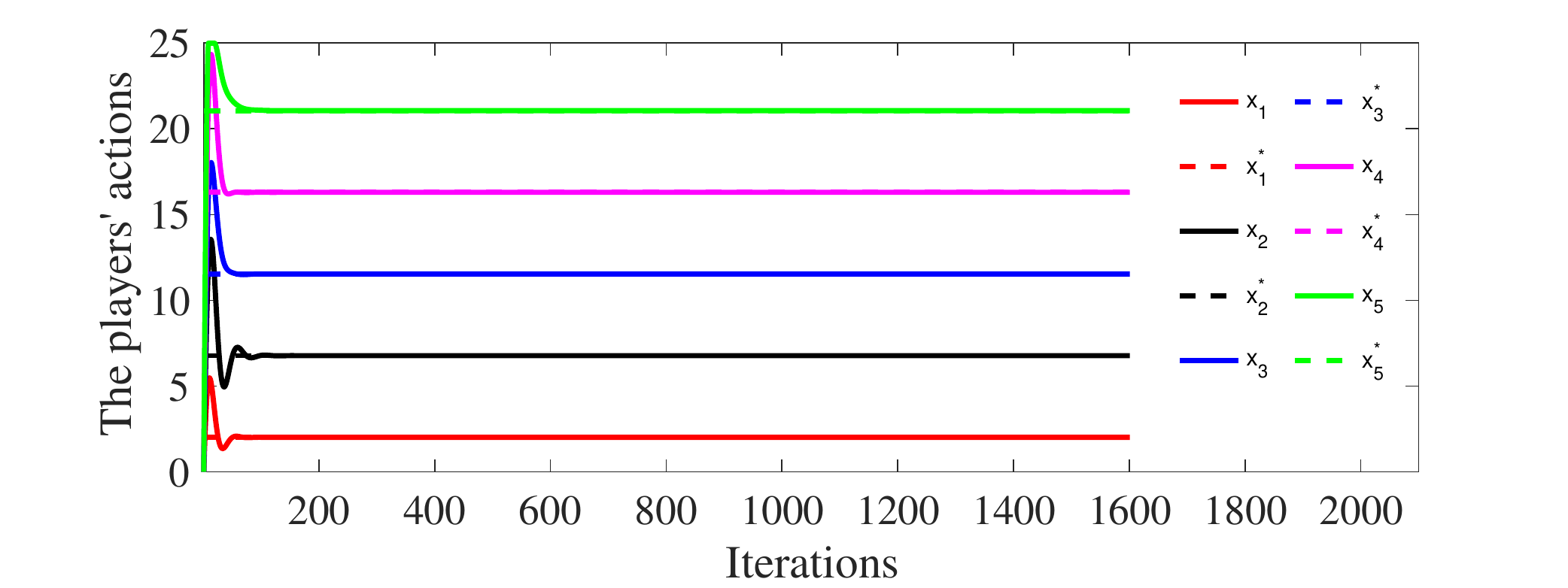} %
\caption{The plot of all players' actions produced by the gradient-based method with constant step-size $\alpha_k= 0.1$.}
\label{fig:constant_stepsize_grad.PNG}
\end{figure}

In this part, we compare the performance of the proposed gradient-free algorithm with its gradient-based counterpart. Specifically, the gradient-based algorithm adopts the same updating laws as in \eqref{eq:update_x} and \eqref{eq:update_y}, but the gradient-free oracle is replaced with the true gradient information. We consider the same problem settings as in section~\ref{sec:simulation}-A under the communication graph as shown in Fig.~\ref{fig:network.png}-(a). For the implementation of the algorithm, the step-size $\alpha_k$ is set to $0.1/\sqrt{k+1}$ and $0.1$, respectively. Figs.~\ref{fig:diminishing_stepsize_rgf.PNG} and \ref{fig:constant_stepsize_rgf.PNG} present the players' actions generated by the proposed gradient-free algorithm with both diminishing step-size and constant step-size, respectively. For the gradient-based counterpart, the convergence results of the players' actions with both diminishing step-size $\alpha_k = 0.1/\sqrt{k+1}$ and constant step-size $\alpha_k = 0.1$ are plotted in Figs.~\ref{fig:diminishing_stepsize_grad.PNG} and \ref{fig:constant_stepsize_grad.PNG}, respectively.

Comparing Figs.~\ref{fig:diminishing_stepsize_rgf.PNG} and \ref{fig:diminishing_stepsize_grad.PNG} for diminishing step-size, and Figs.~\ref{fig:constant_stepsize_rgf.PNG} and \ref{fig:constant_stepsize_grad.PNG} for constant step-size, it can be observed that the convergence speed of the gradient-based algorithm is generally faster than its gradient-free counterpart for both diminishing and constant step-size scenarios. This result is reasonable because the gradient-based algorithm has direct access to the true gradient, where the structure information is included. On the other hand, the faster speed implies more aggressive updates in the process, leading to a relatively larger overshoot, which can be moderated by a smaller step-size.

\section{Conclusions}\label{sec:conclusion}
We have developed a gradient-free distributed Nash equilibrium seeking algorithm for non-cooperative games among a group of players under a directed and strongly connected communication graph. The proposed algorithm does not require the knowledge on the explicit analytical expression of the cost function and allows the problem to be non-smooth. The convergence of the proposed algorithm to the Nash equilibrium has been rigorously studied for both diminishing and constant step-sizes, respectively. Specifically, by choosing a diminishing smoothing parameter, we have shown the convergence to the exact Nash equilibrium for diminishing step-size, and the neighborhood of the Nash equilibrium for constant step-size, in which the gap is proportional to the step-size. Finally, we have illustrated the performance of the algorithm through a numerical example in the application of HVAC system. An outlook to the future research can be the consideration of the dynamical systems in the gradient-free settings. Specifically, instead of having static unknown cost functions, players may follow some dynamics, where the exact model of these dynamics are unknown to the players.









\bibliographystyle{IEEEtran}
\bibliography{d_ne_rgf_reference}

\end{document}